\documentclass[a4paper,11pt]{article}

\usepackage{citesort,geometry,enumerate,amsmath,amssymb,latexsym,theorem,calc,xypic}
{\theorembodyfont{\rmfamily}}
{\theorembodyfont{\rmfamily}}

{\theorembodyfont{\rmfamily}}

\newenvironment{proof}{{\noindent \bf Proof:}}{\hfill $\Box$

\mbox{}}

\def\ge{\geqslant}
\def\bB{\mathbb{B}}
\def\pr{\mathrm{pr}}

\begin{document}

\title{Cubical abelian groups with connections\\ are
equivalent to chain complexes}

\author{Ronald Brown\thanks{email: r.brown@bangor.ac.uk}, \\
 Mathematics Division \\ School of Informatics,  \\ University of
Wales, Bangor \\Gwynedd LL57 1UT,  U.K. \and  Philip
J. Higgins\thanks{email: p.j.higgins@durham.ac.uk}, \\
Department of Mathematical Sciences, \\ Science Laboratories, \\ South Rd., \\
Durham, DH1 3LE,  U.K.}

\maketitle

\begin{center}
University of Wales, Bangor, Maths Preprint  02.24
\end{center}
\begin{abstract}
The theorem of the title  is deduced from  the equivalence between
crossed complexes and cubical $\omega$-groupoids with connections
proved by the authors in 1981. In fact we prove the equivalence of
five categories defined  internally to an additive category with
kernels.
\end{abstract}

\section*{Introduction}

The theorem of the title is shown to be a consequence of the
equivalence between crossed complexes and cubical
$\omega$-groupoids with connections proved by us in \cite{BH3}. We
assume the definitions given in \cite{BH3}. Thus this paper is a
companion to others, for example \cite{T1}, which show that a
deficit of the traditional theory of cubical sets and cubical
groups has been the lack of attention paid to the ``connections'',
defined in \cite{BH3}. Indeed the traditional degeneracies of
cubical theory identify certain opposite faces of a cube, unlike
the degeneracies of simplicial theory which identify adjacent
faces. The connections allow for a fuller analogy with the methods
available for simplicial theory by giving forms of `degeneracies'
which identify adjacent faces of cubes. They are used in
\cite{BH3} and \cite{ABS} to give a definition of a `commutative
cube'.

Part of the interest of these results is that the family of
categories equivalent to that of crossed complexes can be regarded
as a foundation for a non-abelian approach to algebraic topology
and the cohomology of groups. These results show that a form of
abelianisation of these categories leads to well-known structures.

\section*{Crossed complexes internal to an additive category with kernels}

The basic elements of what we say next are well known, but are
given for completeness.

Suppose we are given an action of a group $P$ on the right of a
group $M$ such that the action $\phi : M \times P \to M$ is a
morphism of groups. Then, as is well known, the action is trivial.
The proof is easy: let $m\in M, p\in P.$ Then $m^p = \phi (m,p) =
\phi (m,1) \phi (1,p) = m^1 1^p = m. $ It follows that a crossed
module internal to the category of groups is just a morphism of
abelian groups.

We need to consider below the more general case of crossed modules
over groupoids. Internally to the category of groups, these are
more complicated; but internally to the category of abelian groups
they are again equivalent to morphisms of abelian groups. This
result is essentially in \cite{Gr1}.

\vspace{2ex}

\noindent{\bf Theorem} {\it Let $\cal A$ be an additive category
with kernels. The following categories, defined internally to
$\cal A $, are equivalent.
\begin{enumerate}[$\bB_1$ :] \item  The category of chain complexes. \item
The category of crossed complexes \item  The category of cubical
sets with connections. \item The category of cubical
$\omega$-groupoids with connections. \item  The category of
globular $\omega$-groupoids.
\end{enumerate}}

 \begin{proof} By working on the morphism sets, we can as usual
assume that we are working in the category of abelian groups. Note
that the theorem of the title follows from the equivalence $\bB_3
\simeq \bB_1$.

 $\bB_1 \simeq \bB_2:$ By a chain complex we shall always
mean a sequence of objects and morphisms $\delta : A_n \to
A_{n-1}, \; n \ge 1,$ such that $\delta \delta = 0.$ Let $C$ be a
crossed complex internal to $\cal A$. The associated chain complex
$\alpha C$
will be defined by  \begin{align*} (\alpha C)_0& =  C_0, \\
(\alpha C)_1 & =  \mbox{Ker}\; (\delta_0 : C_1 \to C_0), \\
(\alpha C)_n & =  C_n(0), \; n \ge 2. \\ \intertext{  The crossed
complex $\beta A$   associated to a chain complex $A$
will be defined by} (\beta A)_0 & =  A_0, \\
(\beta A)_1 & =  A_0 \times A_1, \\ (\beta A)_n & =  A_0 \times
A_n, \; n \ge 2 . \end{align*} The groupoid structure on $\beta A$
in dimension 1 is defined as usual by $\delta_0 = \pr_1, \;
\delta_1 = \pr_1 +(\partial \circ \pr_2), \; $ and with
composition $(a,b) + (a+\partial b,c) = (a, b+c). $ The structure
on  $(\beta A)_n $ for $n \ge 2$ is that the only addition is
$(a,b) + (a,c) = (a,b+c). $ The operation of $(\beta A)_1$ on
$(\beta A)_n, \; n \ge 2, $ is $(a,b)^{(a,c)} = (a +
\partial c , b). $ This gives our first equivalence, between
chain complexes and crossed complexes.

 $\bB_2 \simeq \bB_3:$ An equivalence between crossed
complexes and cubical $\omega$-groupoids with connections
internally to the category of sets is established in \cite{BH3}.
Although choices are involved in this, the end result is a natural
equivalence. It follows that this can be applied internally to a
category $\cal A$, simply by applying it to the morphism sets
${\cal A}(X, A)$ for all objects $X$ of $\cal A.$ This yields our
equivalence between crossed complexes and cubical
$\omega$-groupoids with connections internal to $\cal A.$

 $\bB_2 \simeq \bB_5:$ This follows, in a similar way,
from the equivalence between crossed complexes and globular
$\omega$-groupoids proved in \cite{BH2}. (Reference \cite{Bo} is
relevant to the equivalence $\bB_1 \simeq \bB_5.$

 $\bB_3 \simeq \bB_4:$ Let $K$ be a cubical abelian group
with connections, in the sense of \cite{BH3}.

\vspace{2ex}

\noindent {\bf Lemma} {\em If $G$ is an abelian group, and if $s,t
: G \to G$ are endomorphism of $G$ such that $st=s, \; ts=t, $
then we can define a groupoid structure on $G$ with source and
target maps $s,t$ by
\begin{equation*}\label{comp} g \circ h = g - tg + h,
\end{equation*}for $g,h \in G$ with $tg=sh$, and this defines on
$G$ the structure of groupoid internal to abelian groups. }

This result comes from  \cite{Gr1}, and is also a special case of
a non-abelian result on cat$^1$-groups \cite{Lo}, where the
condition $[ \mbox{Ker} \; s, \; \mbox{Ker} \; t] = 1$ is
required, and is here trivially satisfied. This result can be
applied to $K_n, \; n \ge 1, $ and for each $i = 1,\ldots ,n,$
with $s_i= \epsilon_i
\partial^0_i, \; t_i= \epsilon_i \partial^1_i, $ giving $n$ compositions
and so a cubical complex with compositions and connections in the
sense of \cite{ABS,BH3}. The interchange law is easily verified,
and there remains essentially only the transport law for the
connections, which is again simple, showing that $K$ is now a
cubical $\omega$-groupoid with connections. It is easy to see that
the  functor thus defined is adjoint to the forgetful functor
$\bB_4 \to \bB_3$.
\end{proof}


\begin{thebibliography}{99}
\bibitem{ABS}  {\sc Al-Agl, F.A., Brown, R. and Steiner, R.}, `Multiple categories: the
equivalence between a globular and cubical approach', {\em
Advances in Mathematics}    170 (2002) 71-118.
\bibitem{Bo}{\sc Bourn, D.},
`Another denormalization theorem for the abelian chain complexes',
{\em J. Pure Appl. Algebra} 66 (1990) 229-249.
 \bibitem{BH2} {\sc
Brown, R. and  Higgins, P.J. },  `The  equivalence  of
$\infty$-groupoids   and   crossed complexes', {\em Cah. Top.
G\'eom. Diff.} 22 (1981) 371-386. \bibitem{BH3} {\sc  Brown, R.
and Higgins, P.J.}, `The algebra of cubes', {\em J. Pure  Appl.
Algebra}  21 (1981)  233-260.   \bibitem{Gr1} {\sc Grothendieck,
A.}, `Cat\'egories cofibr\'ees additives et complexe cotangent
              relatif', {\em Springer Lecture Notes in Math.} 79
              (1968) Springer-Verlag, Berlin, 167pp.
\bibitem{Lo} {\sc   Loday, J.-L.},  `Spaces with finitely
many non-trivial homotopy groups', {\em J. Pure Appl. Algebra} 24
(1982) 179-202.

\bibitem{T1} {\sc Tonks, A.P.}, `Cubical groups which are Kan',
{\em J. Pure Appl. Algebra} 81 (1992) 83-87.
\end{thebibliography}
\end{document}